\documentclass[11pt,twoside]{article}
\usepackage{amssymb}
\usepackage{indentfirst}
\usepackage{dsfont}
\usepackage{ntheorem}
\usepackage{amsmath}

\begin{document}
\baselineskip 14.2pt \noindent \thispagestyle{empty}

\begin{center} \large{\bf Rigidity of Proper Holomorphic Self-mappings of
  the Pentablock}
\end{center}

\begin{center}
\noindent\text{Guicong Su,\; Zhenhan Tu$^{*}$\; and \; Lei Wang}\\
\vskip 4pt \noindent\small {School of Mathematics and Statistics,
Wuhan University, \\
Wuhan, Hubei 430072, P.R. China} \\
\noindent\text{Email: suguicong@whu.edu.cn (G. Su),\;
zhhtu.math@whu.edu.cn (Z. Tu)} \\  and\; wanglei2012@whu.edu.cn (L.
Wang)
\renewcommand{\thefootnote}{{}}
\footnote{\hskip -16pt {$^{*}$Corresponding author. \\ } }
\end{center}

\begin{center}
\begin{minipage}{13cm}
{\bf Abstract.}
 {\small The pentablock is a Hartogs domain in
$\mathbb{C}^{3}$ over the symmetrized bidisc in $\mathbb{C}^{2}$.
The domain is a bounded inhomogeneous pseudoconvex domain, and does
not have a $\mathcal{C}^{1}$ boundary. Recently, Agler-Lykova-Young
constructed a special subgroup of the group of holomorphic
automorphisms of the pentablock, and Kosi\'{n}ski completely
described the group of holomorphic automorphisms of the pentablock.
The purpose of this paper is to prove that any proper holomorphic
self-mapping of the pentablock must be an automorphism. \vskip 5pt
 {\bf Key words:} Automorphisms, Hartogs domains,  proper holomorphic self-mappings, symmetrized bidisc, the pentablock. \vskip 5pt
 {\bf 2010 Mathematics Subject Classification:} Primary 32A07,  32H35.}
\end{minipage}
\end{center}

\newtheorem{rem}{Remark}
\newtheorem{lem}{Lemma}
\newtheorem{thm}{Theorem}
\newtheorem{pro}{Proposition}
\newcommand{\trace}{\mathop{\mathrm{tr}}}
\section{Introduction}

Let $\mathbb{C}^{2\times 2}$ denote the space of $2\times 2$ complex
matrices, with the usual operator norm, i.e., for a matric $A \in
\mathbb{C}^{2\times 2}$,
$$||A||:=\sup\left\{ ||zA||/||z||:\; z\in \mathbb{C}^{2},
z\not=0\right\}, $$ in which $\mathbb{C}^{2}$ is equipped with the
standard Hermitian norm. Recently Agler-Lykova-Young \cite{Aglerz}
introduced the bounded domain $\mathcal{P}$ by
$$\mathcal{P}:=\left\{(a_{21}, \trace{A}, \det{A}): A=[a_{ij}]_{i, j=1}^{2}\in\mathbb{B} \right\},$$
where
$$\mathbb{B}:= \left\{A\in \mathbb{C}^{2\times 2}:\; ||A||<1 \right\} $$
denotes the open unit ball in the space $\mathbb{C}^{2\times 2}$
with the usual operator norm. So $\mathcal{P}$ is an image of
$\mathbb{B}$ under the holomorphic mapping $A=[a_{ij}]\mapsto
(a_{21}, \trace{A}, \det{A}).$ The domain $\mathcal{P}$ is called
the {\it pentablock} in Agler-Lykova-Young  \cite{Aglerz}, as
$\mathcal{P}\cap\mathbb{R}^{3}$ is a convex body bounded by five
faces, in which three of them are flat and two are curved.

The pentablock $\mathcal{P}$ is polynomially convex and starlike
about the origin, but neither circled nor convex, and  does not have
a $\mathcal{C}^{1}$ boundary (see Agler-Lykova-Young \cite{Aglerz}).
The pentablock is a bounded inhomogeneous domain (see Th. 15 in
Kosi\'{n}ski \cite{Kosin}). For the complex geometry and function
theory of the pentablock $\mathcal{P}$, see Agler-Lykova-Young
\cite{Aglerz} and Kosi\'{n}ski \cite{Kosin} for details.

The pentablock $\mathcal{P}$ arises in connection with the {\it
structured singular value}, a cost function on matrices introduced
by control engineers in the context of robust stabilization with
respect to modelling uncertainty (e.g., see Doyle \cite{Doyle}). The
structured singular value is denoted by $\mu$, and engineers have
proposed an interpolation problem called the $\mu$-{\it synthesis
problem} that arises from this source. Attempts to solve cases of
this interpolation problem have also led to the study of two other
domains,  the {\it symmetrised disc} (e.g., see \cite{Agler1, Edig})
and the {\it tetrablock} (e.g., see \cite{Edig0, Y}),
 in $\mathbb{C}^{2}$ and $\mathbb{C}^{3}$ respectively,
 which have turned out to have many properties of interest to specialists in several complex variables (e.g., see \cite{Agler2, Edig,
Kosi}) and to operator theorists (e.g., see \cite{Bha, Sar}).

Here and throughout the paper, $\mathbb{D}$ denotes the open unit
disc in the complex plane, additionally by $\mathbb{T}$ we shall
denote the unit circle.

The pentablock is closely related to the symmetrized bidisc
$\mathbb{G}_{2}$, which is a bounded domain in $\mathbb{C}^{2}$ as
follows:
$$\mathbb{G}_{2}:=\left\{(\lambda_{1}+\lambda_{2}, \lambda_{1}\lambda_{2})\in\mathbb{C}^{2}: \;\lambda_{1},
 \lambda_{2}\in\mathbb{D}\right\}.$$
For $(s, p)\in\mathbb{C}^{2}$, it is easy to check (see Th. 2.1 in
\cite{Aglerz}) that $(s, p)\in\mathbb{G}_{2}$ if and only if
$|s-\bar{s}p|+|p|^{2}<1$. So the symmetrized bidisc $\mathbb{G}_{2}$
can be also described as $$\mathbb{G}_{2}=\left\{(s,
p)\in\mathbb{C}^{2}:\; |s-\bar{s}p|+|p|^{2}<1\right\}.$$ The
symmetrized bidisc is important since it is the first known example
of a bounded pseudoconvex domain for which the Carath\'eodory and
Lempert functions coincide, but which cannot be exhausted by domains
biholomorphic to convex ones (see Costara \cite{C} and Edigarian
\cite{Edig1}).

Edigarian-Zwonek \cite{Edig} gave the characterization of proper
holomorphic self-mappings of the symmetrized polydisc, which, in the
special case of the symmetrized bidisc, reduces the result  as
follows.

\begin{thm}(Edigarian-Zwonek \cite{Edig})\label{33} Let $f:
\mathbb{G}_{2}\rightarrow \mathbb{G}_{2}$ be a holomorphic mapping.
Then $f$ is proper if and only if there exists a finite Blaschke
product $B$ such that
$$f(\lambda_{1}+\lambda_{2}, \lambda_{1}\lambda_{2})= \left(
B(\lambda_{1})+B(\lambda_{2}), B(\lambda_{1})B(\lambda_{2}))\;\;
(\lambda_{1},
 \lambda_{2}\in\mathbb{D}   \right).
$$
In particular, $f$ is an automorphism if and only if
$$f(\lambda_{1}+\lambda_{2}, \lambda_{1}\lambda_{2})= \left(
\nu(\lambda_{1})+\nu(\lambda_{2}),
\nu(\lambda_{1})\nu(\lambda_{2}))\;\;  (\lambda_{1},
 \lambda_{2}\in\mathbb{D}   \right),
$$
where $\nu$ is an automorphism of the unit disc $\mathbb{D}$. As
$\mbox{\rm Aut}(\mathbb{G}_{2})$ does not act transitively on
$\mathbb{G}_{2}$, the symmetrized bidisc is inhomogeneous.
\end{thm}

By definition of the domain $\mathcal{P}$,  the pentablock is a
Hartogs domain in $\mathbb{C}^{3}$ over the symmetrized bidisc
$\mathbb{G}_{2}$. Indeed, it is clear from the definition that
$\mathcal{P}$ is fibred over $\mathbb{G}_{2}$ by the map
$$(a, s, p) \mapsto (s, p),$$
since if $A\in \mathbb{B}(\subset \mathbb{C}^{2\times 2})$, then the
eigenvalues of $A$ lie in $\mathbb{D}$ and so $(\trace A, \det A)
\in \mathbb{G}_{2}.$ More precisely (e.g., see Th. 1.1 in
Agler-Lykova-Young \cite{Aglerz}), one has
\begin{eqnarray}\label{00}
\mathcal{P}=\left\{(a, s, p)\in\mathbb{D}\times\mathbb{G}_{2}:
|a|^{2}<{e}^{-u(s, p)}\right\},
\end{eqnarray} where
$$u(s, p)=-2\log\left|1-\frac{\frac{1}{2}s\bar{\beta}}{1+\sqrt{1-|\beta|^{2}}}\right|,$$
in which $\beta=\frac{s-\bar{s}p}{1-|p|^{2}}$. Note the $\exp$ in
display \eqref{00} is natural because it was important in
Kosi\'{n}ski \cite{Kosin} due to Kiselman's results on
$\mathbb{C}$-convex Hartogs domains.

Note we have $u(s,0)=-2\left( \log\left(  1+(1-s\bar
s)^{1/2}\right)-\log 2\right)$ and so
\begin{eqnarray}\label{0000}
\frac{\partial^2 u}{\partial  s \partial\bar s} (0,0)=1/2.
\end{eqnarray}

Theorem 1.1 in Agler-Lykova-Young \cite{Aglerz} also proved
$${e}^{-u(\lambda_{1}+\lambda_{2},\lambda_{1}\lambda_{2})/2}=\frac{1}{2}|1-\lambda_{1}\overline{\lambda_{2}}|+\frac{1}{2}
(1-|\lambda_{1}|^2)^{\frac{1}{2}}(1-|\lambda_{2}|^2)^{\frac{1}{2}}\;\;(
\lambda_{1}, \lambda_{2}\in\mathbb{D})$$ and thus gives another
description of the pentablock as follows:
\begin{eqnarray}\label{01}
\mathcal{P}=\left\{(a,
\lambda_{1}+\lambda_{2},\lambda_{1}\lambda_{2}): \; |a|<
\frac{1}{2}|1-\lambda_{1}\overline{\lambda_{2}}|
+\frac{1}{2}(1-|\lambda_{1}|^2)^{\frac{1}{2}}(1-|\lambda_{2}|^2)^{\frac{1}{2}}
, \lambda_{1},
 \lambda_{2}\in\mathbb{D} \right\}.\end{eqnarray}

In 2014, Agler-Lykova-Young \cite{Aglerz} constructed a special
subgroup of the group of holomorphic automorphisms of the pentablock
as follows.

\begin{thm}\label{AA} (Th. 1.2 in \cite{Aglerz})
All mappings of the form
\begin{eqnarray}\label{1}
\begin{aligned}
f_{\omega,\nu}(a, \lambda_{1}+\lambda_{2}, \lambda_{1}\lambda_{2})=
\left(\frac{\omega(1-|\alpha|^{2})a}{1-\bar\alpha(\lambda_{1}+\lambda_{2})+\bar\alpha^{2}\lambda_{1}\lambda_{2}},
\nu(\lambda_{1})+\nu(\lambda_{2}),
\nu(\lambda_{1})\nu(\lambda_{2})\right),
\end{aligned}
\end{eqnarray}
where $(a, \lambda_{1}+\lambda_{2},
\lambda_{1}\lambda_{2})\in\mathcal{P},$  $\lambda_{1},
\lambda_{2}\in\mathbb{D},$ form a subgroup of the group $\mbox{\rm
Aut}(\mathcal{P})$  of holomorphic automorphisms of the pentablock,
where $\nu$ is a $M\ddot{o}bius$ function of the form
$\nu(\lambda)=\eta\frac{\lambda-\alpha}{1-\bar{\alpha}\lambda}$, in
which $\omega, \eta\in\mathbb{T}$, $\alpha\in\mathbb{D}$.
\end{thm}

Furthermore, Kosi\'{n}ski \cite{Kosin} completely described the
group of holomorphic automorphisms of the pentablock as follows.

\begin{thm}\label{AB} (Th. 15 in \cite{Kosin}) All automorphisms of the form \eqref{1}
comprise the whole group $\mbox{\rm Aut}(\mathcal{P})$ of
holomorphic automorphisms of the pentablock.
\end{thm}

In this paper we study proper holomorphic self-mappings of the
pentablock and  prove that any proper holomorphic self-mapping of
the pentablock must be a mapping of the form \eqref{1} as follows:

\begin{thm}\label{2}
Any proper holomorphic self-mapping of the pentablock must be an
automorphism of the form \eqref{1}.
\end{thm}

This paper will use Theorem 2 to give a unified proof for Theorem 3
and Theorem 4. Generally speaking, a proper holomorphic mapping
between two bounded domains may lead naturally to the geometric
study of a mapping between their boundaries. These researches are
often heavily based on analytic techniques about the mapping on
boundaries (e.g., see Huang \cite{H} for references). As we know,
the pentablock is a bounded inhomogeneous pseudoconvex domain, and
does not have a $\mathcal{C}^{1}$ boundary. The lack of boundary
regularity usually presents a serious analytical difficulty (e.g.,
see Mok-Ng-Tu \cite{MNT}, Tu \cite{Tu1, Tu2} and Tu-Wang \cite{TW}).
The crucial tools used in deducing Theorem \ref{2} involves
holomorphic extendability of proper holomorphic mapping between
quasi-balanced domains whose Minkowsi functions are continuous (see
Kosi\'{n}ski \cite{Kosi}), the description of proper holomorphic
self-mappings of the symmetrized bidisc (see Edigarian-Zwonek
\cite{Edig}), the complex geometry of the boundary of the pentablock
(see Agler-Lykova-Young \cite{Aglerz} and Kosi\'{n}ski
\cite{Kosin}), and the rigidity of proper holomorphic self-mappings
and the description of automorphisms of the ellipsoids (see
Dini-Primicerio \cite{Dini}).

\section{Preliminaries}
As the pentablock $\mathcal{P}$ is a Hartogs domain in
$\mathbb{C}^{3}$ over a symmetrized bidisc $\mathbb{G}_{2}$, it is
important to learn the basic complex geometry of the symmetrized
bidisc $\mathbb{G}_{2}$.

\begin{lem}(Prop. 3.2 in Agler-Lykova-Young \cite{Agler2})\label{3} Let the symmetrized bidisc $\mathbb{G}_{2}$ be defined as
$$\mathbb{G}_{2}:=\left \{(\lambda_{1}+\lambda_{2}, \lambda_{1}\lambda_{2})\in\mathbb{C}^{2}:\;  \lambda_{1},
 \lambda_{2}\in\mathbb{D} \right\}.$$
Then we have the following results.

$(i)$ For $(s, p)\in\mathbb{C}^{2}$, then $(s, p)\in\mathbb{G}_{2}$
if and only if $|s-\bar{s}p|+|p|^{2}<1$;

$(ii)$ For $(s, p)\in\mathbb{C}^{2}$, then $(s,
p)\in\partial\mathbb{G}_{2}$ if and only if $|s|\leq 2$ and
$|s-\bar{s}p|+|p|^{2}=1$.

$(iii)$  By $\partial_{s}\mathbb{G}_{2}$ denote its Shilov boundary.
Then
 $$\partial_{s}\mathbb{G}_{2}=\left\{(\lambda_{1}+\lambda_{2},
 \lambda_{1}\lambda_{2})\in\mathbb{C}^{2}: \; \lambda_{1}, \lambda_{2}\in\mathbb{T} \right\}.$$
\end{lem}

The royal variety $\Sigma$  of the symmetrized bidisc
$\mathbb{G}_{2}$ plays an important role in the study of the
symmetrized bidisc. Recall the royal variety
 $$\Sigma:=\left\{(2\lambda, \lambda^{2})\in\mathbb{C}^{2}:\;
 \lambda\in{\overline{\mathbb{D}}}\right\}\; (\subset \overline{\mathbb{G}}_{2}).$$
Let the mapping $\sigma: \mathbb{D}^{2}\rightarrow\mathbb{G}_{2}$ be
defined by
$$\sigma(\lambda_{1}, \lambda_{2}):=(\lambda_{1}+\lambda_{2},
\lambda_{1}\lambda_{2}).$$ Thus the mapping $\sigma:
\mathbb{D}^{2}\rightarrow\mathbb{G}_{2}$ is a proper holomorphic
mapping. Note $\sigma$ is well-defined on $\mathbb{C}^{2},$ and
$$\sigma:\; \mathbb{C}^{2} \setminus \left\{(\lambda, \lambda):\;
 \lambda\in\mathbb{C}\right\}
\longrightarrow\mathbb{C}^{2}\setminus  \left\{(2\lambda,
\lambda^2):\;
 \lambda\in\mathbb{C}\right\}$$ is a holomorphic covering. Thus we have
$\sigma(\partial \mathbb{D}^{2})= \partial \mathbb{G}_{2},$ and the
boundary part
$\partial\mathbb{G}_{2}\setminus\partial_{s}\mathbb{G}_{2}$ of the
symmetrized bidisc  is a Levi flat part of the boundary.

Now return to the pentablock, by Agler-Lykova-Young \cite{Aglerz},
we have that $\mathcal{P}$ is a domain of holomorphy and
$\mathcal{P}$ does not have a $\mathcal{C}^{1}$ boundary.
Consequently much of the theory of pseudoconvex domains does not
apply to $\mathcal{P}$. But follwing Kosi\'{n}ski \cite{Kosin}, we
obtain some useful boundary properties.

\begin{lem}(Lemmas 11 and 13 in Kosi\'{n}ski \cite{Kosin})  \label{4}
$(i)$ Any point $x$ of
$$\partial_{1}\mathcal{P}: =\left\{(a, s,
p)\in{\mathbb{C}}\times\mathbb{G}_{2}:\;  |a|^{2}=e^{-u(s,
p)}\right\} \; (=
\partial\mathcal{P} \cap ({\mathbb{C}}\times\mathbb{G}_{2}))
$$
is a smooth point of the boundary $\partial\mathcal{P}.$
Moreover the rank of the Levi form of a defining function of
$\partial\mathcal{P}$ at the point $x$ restricted to the complex
tangent space is equal to $1$; As $u$ is not a pluriharmonic
function, $\partial_{1}\mathcal{P}$ is not Levi flat.

$(ii)$ The boundary part
$$\partial_{2}\mathcal{P}: =\left\{(a, s,
p)\in\mathbb{D}\times  \partial\mathbb{G}_{2}: \; (s,
p)\in\partial\mathbb{G}_{2}\setminus\partial_{s}\mathbb{G}_{2},\;
|a|^{2}<e^{-u(s, p)}\right\}
$$ is a Levi flat part of the boundary
$\partial\mathcal{P}$ and $\partial_{2}\mathcal{P}\subset
\partial\mathcal{P} \cap ({\mathbb{C}}\times (\partial\mathbb{G}_{2}\setminus\partial_{s}\mathbb{G}_{2}))
$.
\end{lem}

Note that we have ${e}^{-u(s, p)}>0$ on
$\partial{\mathbb{G}}_{2}\setminus \Sigma$ by combining \eqref{00}
and \eqref{01}. Therefore, ${e}^{-u(s, p)}>0$ on
$\partial\mathbb{G}_{2}\setminus
\partial_{s}\mathbb{G}_{2}.$ So
$\partial_{2}\mathcal{P}$ is obviously
 a Levi flat part of the boundary
$\partial\mathcal{P}$, as
$\partial\mathbb{G}_{2}\setminus\partial_{s}\mathbb{G}_{2}$ is a
Levi flat part of the boundary $\partial\mathbb{G}_{2}$.

% \begin{lem}(Th. 8.4 in Agler-Lykova-Young \cite{Aglerz})\label{5} $$\partial_{s}\mathcal{P}: =\left\{(a, s, p)\in\mathbb{C}^{3}:\;  (s,
%  p)\in\partial_{s}\mathbb{G}_{2},\;  |a|^{2}=e^{-u(s, p)}\right\}$$
% is the Shilov boundary of the pentablock $\mathcal{P}$. \end{lem}

We will use the notion of quasi-circular domains. Let $m_1,\cdots,
m_n$ be relatively prime natural numbers. Recall that a domain
$D\subset\mathbb{C}^{n}$ is said to be $(m_1,\cdots,m_n)$-circular
(shortly quasi-circular) if
$$(\lambda^{m_1}z_1, \cdots, \lambda^{m_n}z_n)\in D$$
for all $\lambda\in \mathbb{T}$ and $z=(z_1,\cdots, z_n)\in D.$ If
the relation holds for all $\lambda\in\mathbb{D}$, then $D$ is said
to be $(m_1,\cdots,m_n)$-balanced (shortly quasi-balanced). After a
minor modification of the argument in Bell \cite{Bell}, Kosi\'{n}ski
\cite{Kosi} get the holomorphic extendability of the proper
holomorphic mappings as follows.

\begin{lem}(Lemma 6 in Kosi\'{n}ski \cite{Kosi})\label{KK}
Let $D,\; G$ be bounded domains in $\mathbb{C}^{n}.$ Suppose that
$G$ is $(m_1,\cdots,m_n)$-circular and contains the origin. Assume
moreover that the Bergman kernel function $K_{D}(z,
\overline{\xi})\;(z,\xi\in D)$ associated to $D$ satisfies the
following property: for any open, relatively compact subset $E$ of
$D$ there is an open set $U=U(E)$ containing $\overline D$ such that
$K_{D}(z, \overline{\xi})$ extends to be holomorphic on $U$ as a
function of $z$ for each ${\xi} \in E$. Then any proper holomorphic
mapping $f:\; D\rightarrow G$ extends holomorphically to a
neighborhood of $\overline D$.
\end{lem}

Now we consider the proper holomorphic self-mappings of the
pentablock $\mathcal{P}$. It follows from Agler-Lykova-Young
\cite{Aglerz} that the pentablock $\mathcal{P}$ is
$(1,1,2)$-balanced and $\mathcal{P}_r=\{(ra, rs, r^2 p):\ (a,s,p)\in
\mathcal{P}\}$ are relatively compact in $\mathcal{P}$ for any
$0<r<1$. Thus by Remark 7 in Kosi\'{n}ski \cite{Kosi}, (quasi)
Minkowski functional of the pentablock $\mathcal{P}$ is continuous,
and thus the pentablock $\mathcal{P}$ fulfils the assumptions of $D$
in Lemma \ref{KK}. Therefore, Lemma \ref{KK} implies the holomorphic
extendability of the proper holomorphic self-mappings of
$\mathcal{P}$ as follows.

\begin{lem}\label{6}
Any proper holomorphic self-mapping of the pentablock $\mathcal{P}$
extends holomorphically to a neighborhood of
$\overline{\mathcal{P}}$.
\end{lem}

Since $\mathcal{P}$ is a Hartogs domain over a symmetrized bidisc,
we define
$$\Omega: = \left\{(a, 2\lambda, \lambda^{2})\in{\mathcal{P}}: \; \lambda\in\mathbb{D} \right\},$$ i.e.,
$\Omega=\left\{(a, 2\lambda, \lambda^{2}):\; |a|<1-|\lambda|^{2},\;
\lambda\in\mathbb{D}\right\}$, and set $\widetilde{\Omega}:
=\left\{(a, s):\; (a, s, p)\in\Omega\right\}.$ Then we have
$$\widetilde{\Omega}=\left\{(a,
s)\in \mathbb{C}^{2} : |a|+\frac{|s|^{2}}{4}<1\right\}$$ is an
ellipsoid. By the rigidity of proper holomorphic self-mappings and
the description of automorphisms of the ellipsoids (see Corollary
1.2 in Dini-Primicerio \cite{Dini}), we get the following lemma as
follows.

\begin{lem}(Dini-Primicerio \cite{Dini})\label{7} Let $\varphi : \widetilde{\Omega}\rightarrow\widetilde{\Omega}$ be a proper holomorphic mapping.
Then $\varphi$ is an automorphism of $\widetilde{\Omega}$. Moreover,
if $\varphi(0, 0)= (\xi, 0)$ for some $|\xi|<1$, then there exists
$\lambda_1, \lambda_2 \in\mathbb{T}$ such that $\varphi(a ,s)\equiv
(\lambda_1 a, \lambda_2 s)$
 on $\widetilde{\Omega}$.
\end{lem}

\section{Proof of Theorem \ref{2}}
Let $f=(f_{1}, f_{2}, f_{3}): \mathcal{P}\rightarrow\mathcal{P}$ be
a proper holomorphic mapping. By Lemma \ref{6}, $f$ extends
holomorphically to a neighborhood $V$ of $\overline{\mathcal{P}}$
with $f(\partial\mathcal{P})\subset\partial\mathcal{P}$. Define
$$ S:=\{\xi\in V: J_{f}(\xi)=0\},$$
where $J_{f}(\xi)=\det({\frac{\partial f_{i}}{\partial \xi_{j}}})$
is the complex Jacobian determinant of $f(\xi)$ $(\xi\in V).$

Step 1. Consider the mapping $\Psi: \mathbb{G}_{2}\rightarrow
\mathbb{G}_{2}$ defined by
\begin{equation}\label{8}
\Psi: (s,p)\mapsto{(f_{2}(0, s, p), f_{3}(0, s, p))}.
\end{equation}
Then $\Psi$ extends holomorphically to the closure
$\overline{\mathbb{G}}_{2}$. We will prove that $\Psi:
\mathbb{G}_{2}\rightarrow \mathbb{G}_{2}$ is a proper holomorphic
self-mapping of $\mathbb{G}_{2}$ here.

Since $\partial_{2}\mathcal{P}$ is a Levi flat part of the boundary
of the pentablock by Lemma \ref{4}, by the local biholomorphism of
$f$ on $\partial_{2}\mathcal{P}\setminus S$, we have that, for any
$(a, s, p)\in
\partial_{2}\mathcal{P}\setminus S$,  $f(a, s, p)$ lies in a Levi
flat part of the boundary of the pentablock and thus  $f(a, s, p)\in
\partial\mathcal{P}\setminus \partial_{1}\mathcal{P}\subset \mathbb{C}\times \partial\mathbb{G}_2$.

By Lemma \ref{3}, we have
$$\partial\mathbb{G}_{2}= \{(s, p):
|s|\leq 2 \; \mbox{and} \;|s-\bar{s}p|+|p|^{2}=1\}.$$ Thus, for $(a,
s,p)\in\partial_{2}\mathcal{P}\setminus S$, we have
$$|f_{2}(a, s, p)-\overline{f_{2}(a, s, p)}f_{3}(a, s, p)| +|f_{3}(a,
s, p)|^{2}=1.$$
Because of the density of
$\partial_{2}\mathcal{P}\setminus S$ in $\partial_{2}\mathcal{P}$
and the continuity of $f$ on $\overline{\mathcal{P}}$, we conclude
\begin{equation}\label{09}
|f_{2}(a, s, p)-\overline{f_{2}(a, s, p)}f_{3}(a, s, p)| +|f_{3}(a,
s, p)|^{2}=1
\end{equation}
for all $(a, s,p)\in\partial_2\mathcal{P}.$

Note that $(0, s,p)\in\partial_2\mathcal{P}$ for all
$(s,p)\in\partial\mathbb{G}_{2}\setminus \partial_s\mathbb{G}_{2},$
as ${e}^{-u(s, p)}>0$ on $\partial\mathbb{G}_{2}\setminus
\partial_s\mathbb{G}_{2}$.
This means
$$|f_{2}(0, s, p)-\overline{f_{2}(0, s, p)}f_{3}(0, s, p)| +|f_{3}(0,
s, p)|^{2}=1$$ for all $(s,p)\in\partial\mathbb{G}_{2}\setminus
\partial_s\mathbb{G}_{2},$ and then holds for all
$(s,p)\in\partial\mathbb{G}_{2}$ by the continuity of $f$ on
$\overline{\mathcal{P}}$. Obviously $|f_{2}(0, s, p)|\leq 2$ for all
$(s,p)\in\partial\mathbb{G}_{2}$. So $\Psi$ maps
$\partial\mathbb{G}_{2}$ into $\partial\mathbb{G}_{2}$. Thus we get
that
$$\Psi:\; (s,p)\mapsto{(f_{2}(0, s, p), f_{3}(0, s, p))}$$
is a proper holomorphic self-mapping of $\mathbb{G}_{2}$. So,
Theorem 1 implies which there is a finite Blaschke product $b$ such
that
\begin{equation}\label{9}
(f_{2}(0, \lambda_{1}+\lambda_{2}, \lambda_{1}\lambda_{2}), f_{3}(0,
\lambda_{1}+\lambda_{2},
\lambda_{1}\lambda_{2})=(b(\lambda_{1})+b(\lambda_{2}),
b(\lambda_{1})b(\lambda_{2}))
\end{equation}
for all $\lambda_{1}, \lambda_{2}\in\mathbb{D}.$

 Step 2. Here we will
prove that the proper holomorphic mapping $f=(f_{1}, f_{2}, f_{3}):
\mathcal{P}\rightarrow\mathcal{P}$ must map a fiber to a fiber over
$\mathbb{G}_{2}$.

Fix $(s, p)\in\partial_{s}\mathbb{G}_{2}$ with $s^{2}\neq 4p$ (i.e.,
$(s, p)\in\partial_{s}\mathbb{G}_{2}\setminus \Sigma$). As
${e}^{-u(s, p)}>0$ on $\partial_s\mathbb{G}_{2}\setminus \Sigma$
(this can be seen by combining \eqref{00} and \eqref{01}), we have
$(a, s, p)\in\partial\mathcal{P}$ for $|a|^{2}<e^{-u(s, p)}$.
However, by combining (\ref{8}) and (\ref{9}), we have
$$|f_{3}(0, s, p)|=1$$
for any $(s, p)\in\partial_{s}\mathbb{G}_{2}$, as
$|b(\lambda_{1})b(\lambda_{2})|=1$ for $\lambda_{1},\lambda_{2}\in
\mathbb{T}$. This means that $f_{3}(a, s, p)$ $(|a|^{2}<e^{-u(s,
p)})$ attains its maximum modulus at the point $a=0$. Thus $f_{3}(a,
s, p)$ is independent of $a$ for fixed $(s,
p)\in\partial_{s}\mathbb{G}_{2}$ with $s^{2}\neq 4p$.

Now still fix $(s, p)\in\partial_{s}\mathbb{G}_{2}$ with $s^{2}\neq
4p$. Then we have $f_{3}(a, s, p)\equiv e^{i\theta}$ for some
$\theta\in\mathds{R}$. Because of the density of
$\partial_{2}\mathcal{P}$ in $\partial\mathcal{P} \cap
({\mathbb{C}}\times
\partial\mathbb{G}_{2})$ and
the continuity of $f$ on $\overline{\mathcal{P}}$, from \eqref{09}
we conclude
$$
|f_{2}(a, s, p)-\overline{f_{2}(a, s, p)}f_{3}(a, s, p)| +|f_{3}(a,
s, p)|^{2}=1
$$
for all $(a, s,p)\in\partial\mathcal{P} \cap ({\mathbb{C}}\times
\partial\mathbb{G}_{2}).$ Then, we have
$f_{2}(a,s,p)=\overline{f_{2}(a,s,p)}e^{i\theta}$  for all $
|a|^{2}<e^{-u(s, p)}$. Thus $f_{2}(a, s, p)$ is independent of $a$
for fixed $(s, p)\in\partial_{s}\mathbb{G}_{2}$ with $s^{2}\neq 4p$
as well.

Therefore, $f_{2}(a, s, p)$  and $f_{3}(a, s, p)$ are independent of
$a$ for $ |a|^{2}<e^{-u(s, p)}$ with fixed $(s,
p)\in\partial_{s}\mathbb{G}_{2}$ ($s^{2}\neq 4p$). This means that
for all positive integers $k$,
$$\frac{\partial^k f_{2}}{\partial a^k} (0, s, p)=0
\;\;\mbox{and}\;\;\frac{\partial^k f_{3}}{\partial a^k} (0, s,
p)=0$$ for $(s, p)\in\partial_{s}\mathbb{G}_{2}$ ($s^{2}\neq 4p$),
and thus holds for all $(s, p)\in\partial_{s}\mathbb{G}_{2}$ by the
holomorphism of $f$ on $\overline{\mathcal{P}}$. Since
$\partial_{s}\mathbb{G}_{2}$ is the Shilov boundary of
$\mathbb{G}_{2}$, we have that for all positive integers $k$,
$$\frac{\partial^k f_{2}}{\partial a^k} (0, s, p)\equiv 0
\;\;\mbox{and}\;\;\frac{\partial^k f_{3}}{\partial a^k} (0, s,
p)\equiv 0\;\;((s, p)\in\mathbb{G}_{2}).$$  That is, $f_{2}(a, s,
p)$ and $f_{3}(a, s, p)$ defined  on $\mathcal{P}$ are independent
of $a$. So, from \eqref{9}, there exists a finite Blaschke product
$b$ such that
\begin{equation}\label{10}
f(a, \lambda_{1}+\lambda_{2}, \lambda_{1}\lambda_{2})=(f_{1}(a,
\lambda_{1}+\lambda_{2}, \lambda_{1}\lambda_{2}),
 b(\lambda_{1})+b(\lambda_{2}), b(\lambda_{1})b(\lambda_{2}))
\end{equation}
for $(a, \lambda_{1}+\lambda_{2}, \lambda_{1}\lambda_{2})\in
\mathcal{P}$, in which $\lambda_{1}, \lambda_{2}\in\mathbb{D}.$

Step 3. Here we will prove that $f$ is an automorphism of the form
(\ref{1}).

By \eqref{10},  we have $$f(a, 2\lambda, \lambda^2)=(f_{1}(a,
2\lambda, \lambda^2), 2 b(\lambda), b(\lambda)^2)$$ for $(a,
2\lambda, \lambda^2)\in \mathcal{P}$, $\lambda\in\mathbb{D}.$  So
$f$ preserves
$$\Omega\;:= \{(a, 2\lambda, \lambda^{2})\in{\mathcal{P}}: \;
\lambda\in\mathbb{D}\}= \{(a, 2\lambda,
\lambda^{2})\in{\mathbb{C}^3}: \;|a|<1-|\lambda|^{2}, \;
\lambda\in\mathbb{D}\},$$ i.e., $f(\Omega)\subset\Omega$. Note that
the group $\mbox{\rm Aut}(\mathbb{G}_{2})$ of the symmetrized bidisc
$\mathbb{G}_{2}$ acts transitively on $\{(2\lambda,
\lambda^{2})\in\mathbb{C}^{2}:\;
 \lambda\in{{\mathbb{D}}}\} (\subset{\mathbb{G}}_{2}).$
Then, we can take some $\varphi_1\in \mathrm{Aut}(\mathcal{P})$ in
the form (\ref{1}) such that
$$\varphi_1\circ f(0, 0, 0)=\varphi_1(f_1(0, 0, 0), 2b(0), b(0)^2)=(r, 0, 0)$$
for some $0\leq r <1.$ Then
$$F(a, s, p):=\varphi_1\circ f(a, s, p): \; \mathcal{P}\rightarrow
\mathcal{P}$$ is a proper holomorphic self-mapping of $\mathcal{P}$
with  $F(0,0,0)=(r,0,0)$ for some $0\leq r <1.$ So, from \eqref{10}
there exists a finite Blaschke product $\tilde{b}$ with
$\tilde{b}(0)=0$ such that
\begin{equation}\label{11}
F(a, \lambda_{1}+\lambda_{2}, \lambda_{1}\lambda_{2})=(F_{1}(a,
\lambda_{1}+\lambda_{2}, \lambda_{1}\lambda_{2}),\;
 \tilde{b}(\lambda_{1})+\tilde{b}(\lambda_{2}),\; \tilde{b}(\lambda_{1})\tilde{b}(\lambda_{2}))
\end{equation}
for $(a, \lambda_{1}+\lambda_{2}, \lambda_{1}\lambda_{2})\in
\mathcal{P}$, in which $\lambda_{1}, \lambda_{2}\in\mathbb{D}.$

Let
$$\widetilde{\Omega}:
=\{(a, 2\lambda):\; (a, 2\lambda,
\lambda^{2})\in{\mathcal{P}}\}=\{(a, s)\in \mathbb{C}^{2} :
|a|+\frac{|s|^{2}}{4}<1\},$$ an ellipsoid. Then the mapping $\Phi:
\widetilde\Omega\rightarrow \widetilde\Omega $ defined by
$$\Phi: (a, s)\mapsto(F_{1}(a, s, {s^{2}}/{4}), \; 2\tilde{b}(s/2)  )$$
is a proper holomorphic self-mapping of $\widetilde\Omega$ with
$\Phi(0,0)=(r,0).$  Thus by Lemma \ref{7}, there exists
$\eta_1,\eta_2\in \mathbb{T}$ such that
$$\Phi(a, s)=(\eta_1 a,\; \eta_2 s),$$
that is,
$$F_{1}(a, s, {s^{2}}/{4})= \eta_1 a,\;\;\tilde{b}(s)=
\eta_2 s.$$ So, from \eqref{11},  we have
$$F(a, s, p)=(F_{1}(a, s, p),\; \eta_2
s,\;  \eta_2^2p),$$ where $F_{1}(a, s, {s^{2}}/{4})= \eta_1 a.$ Note
that $\varphi_2(a, s,p):=(a/{\eta_1}, \; s/{\eta_2},\; p/{\eta_2^2})
\;(\in \mathrm{Aut}(\mathcal{P}))$ is of the form (\ref{1}) such
that
$$
G(a, s, p):= \varphi_2\circ F(a, s, p)=(G_{1}(a, s, p),\;  s,\; p)$$
is a proper holomorphic self-mapping of $\mathcal{P}$ as well with
$G_{1}(a, s, {s^{2}}/{4})=a.$ Thus the complex Jacobian determinant
$$J_{G}(a,0,0)\equiv 1\;\; (|a|\leq  e^{-u(0, 0)}=1).$$
Since $G(a, s, p)$ is holomorphic on $\overline{\mathcal{P}},$ there
exists a neighborhood $U_0$ of $(0,0)$ in the symmetrized bidisc
$\mathbb{G}_{2}$ such that
$$J_{G}(a,s,p)\not=0\;\mbox{everywhere on}\; \overline{\mathcal{P}} \cap ({\mathbb{C}}\times
U_0).$$ Moreover, as $\frac{\partial^2 u}{\partial s
\partial  \bar s} (0,0)=1/2$ by \eqref{0000}, we may assume
\begin{equation}\label{1111}
\frac{\partial^2 u}{\partial  s \partial\bar s} (s,p)\not=0
\end{equation}
everywhere on $U_0$ also.

Thus we get that $G(a, s, p)=(G_{1}(a, s, p),\;  s,\; p)$ restricted
on ${\mathcal{P}} \cap ({\mathbb{C}}\times U_0)$ is a biholomorphic
self-mapping of ${\mathcal{P}} \cap ({\mathbb{C}}\times U_0)$. So,
for fixed $(s,p)\in U_0,$ $G_1(a, s, p)$ on the disc $|a|^2<
e^{-u(s,p)}$ is a $M\ddot{o}bius$ function, i.e.,
\begin{equation}\label{cc}
G_1(a, s, p)=
e^{\sqrt{-1}\theta(s,p)}\frac{a-G_1^{-1}(0,s,p)}{1-\overline{G_1^{-1}(0,s,p)}e^{u(s,p)}a},\;\;\;\;
(a, s, p)\in {\mathcal{P}} \cap ({\mathbb{C}}\times
U_0).\end{equation}

Then
$$H(a, s, p):= e^{-\sqrt{-1}\theta(s,p)}\left({1-\overline{G_1^{-1}(0,s,p)}e^{u(s,p)}a}\right)$$
is holomorphic in ${\mathcal{P}} \cap
({\mathbb{C}}\times U_0)$, and so
$$e^{-\sqrt{-1}\theta(s,p)}\;(=H(0, s, p))\;\; \mbox{and} \;\; e^{-\sqrt{-1}\theta(s,p)}
\overline{G_1^{-1}(0,s,p)}e^{u(s,p)}\;(= -\frac{\partial
H(a,s,p)}{\partial a})$$ are holomorphic on $U_0$. Thus, from
$|e^{-\sqrt{-1}\theta(s,p)}|\equiv 1$ on $U_0$ we get
$$e^{-\sqrt{-1}\theta(s,p)}\equiv e^{-\sqrt{-1}\theta_0}$$ on  $U_0$.
So $$\overline{G_1^{-1}(0,s,p)}e^{u(s,p)}$$ is holomorphic on $U_0$
also. Assume that ${G_1^{-1}(0,s,p)}\not=0$ everywhere on a small
open ball $V_0$ of $U_0$. Then, from the holomorphism of
$\overline{G_1^{-1}(0,s,p)}e^{u(s,p)}$ on $V_0$, we immediately get
${u(s,p)}$ is pluriharmonic  on $V_0$, a contradiction with
\eqref{1111}. Thus we get ${G_1^{-1}(0,s,p)}\equiv 0$ on  $U_0$.
Therefore,  we get $G_1(a, s, p)= e^{\sqrt{-1}\theta_0}a$ on
${\mathcal{P}} \cap ({\mathbb{C}}\times U_0)$ by  \eqref{cc}  and so
$$G(a, s, p)=(e^{\sqrt{-1}\theta_0}a,\;  s,\; p)$$
on ${\mathcal{P}}$ by the identity principle, which means $G$ is an
automorphism of ${\mathcal{P}}$ in the form (\ref{1}). So we
conclude $f(a, s, p)=\varphi_2^{-1}\circ \varphi_1^{-1}\circ G(a, s,
p)$ is an automorphism of ${\mathcal{P}}$ in the form (\ref{1})
also.  The proof of Theorem 4 is complete.

\vskip 10pt

\noindent\textbf{Acknowledgments.}\quad The authors are grateful to
Professors Ngaiming Mok and Xiaojun Huang for helpful comments. In
addition, the authors would like to thank the referees for many
helpful suggestions. The project was supported by the National
Natural Science Foundation of China (No.11271291).

\end{document}